\documentclass[12pt]{amsart}
\usepackage[all]{xy}
\usepackage[parfill]{parskip}
\usepackage{hyperref}
\usepackage{pinlabel}
\usepackage{enumerate}
\usepackage{verbatim}
\usepackage{amsmath, amscd, graphics, latexsym, times, rlepsf}

\usepackage{subcaption}
\oddsidemargin .25in \theoremstyle{plain}
\newtheorem{dummy}{anything}[section]
\newtheorem{theorem}[dummy]{Theorem}

\theoremstyle{definition}

\newtheorem{remark}[dummy]{Remark}

\newcommand{\Z}{\mathbb{Z}}
\newcommand{\R}{\mathbb{R}}

\newcommand{\Sp}{\mathbb{S}}
\newcommand{\D}{\mathbb{D}}
\renewcommand{\t}{\mathbf t}

\renewcommand{\sp}{\mathbf s}

\newcommand{\T}{\mathbb{T}}

\def\S{\Sigma}

\def\l{\lambda}

\begin{document}

\title{Stein and Weinstein structures on  disk cotangent bundles of surfaces}

\author{Burak Ozbagci}

\address{Department of Mathematics, Ko\c{c} University, Istanbul,
Turkey}
\email{bozbagci@ku.edu.tr}

\subjclass[2000]{}
\thanks{}


\begin{abstract}
In \cite{go}, Gompf describes a Stein domain structure on the disk cotangent bundle of any closed surface $S$, by a Legendrian handlebody diagram.  We prove that Gompf's Stein domain is symplectomorphic to the disk cotangent bundle equipped with its canonical symplectic structure and the boundary of this domain is contactomorphic to the unit cotangent bundle of $S$ equipped with its canonical contact structure. As a corollary, we obtain a surgery diagram for the canonical contact structure on the unit cotangent bundle of $S$.




\end{abstract}

\maketitle

\section{Introduction}\label{sec: intro}

Let $S$ be a closed, connected and smooth surface,  which we do not assume to be orientable. The  disk cotangent bundle $\D T^*S$   of $S$ carries the canonical symplectic structure $\omega_{can} = d \l_{can}$, and the
unit cotangent bundle $\partial (\D T^*S )= \Sp T^*S$  carries the canonical contact  structure $\xi_{can} = \ker (\l_{can}|_{\Sp T^*S})$,  where $\l_{can}$ is the Liouville one form on $T^*S$.    In \cite{go}, Gompf showed that $\D T^*S$ admits the structure of a Stein domain, by explicitly exhibiting $\D T^*S$ as a Legendrian handlebody diagram.

Here we prove that Gompf's Stein domain is symplectomorphic to $(\D T^*S, \omega_{can})$ and the boundary  contact $3$-manifold is contactomorphic to $(\Sp T^*S, \xi_{can})$. As a corollary,
we obtain a contact surgery diagram for  $(\Sp T^*S, \xi_{can})$, using a technique described by Ding and Geiges \cite[Theorem 5]{dg}.

\begin{theorem} \label{thm: gene}

\begin{figure}
\centering
\begin{minipage}{.5\textwidth}
  \centering
  \includegraphics[width=.7\linewidth]{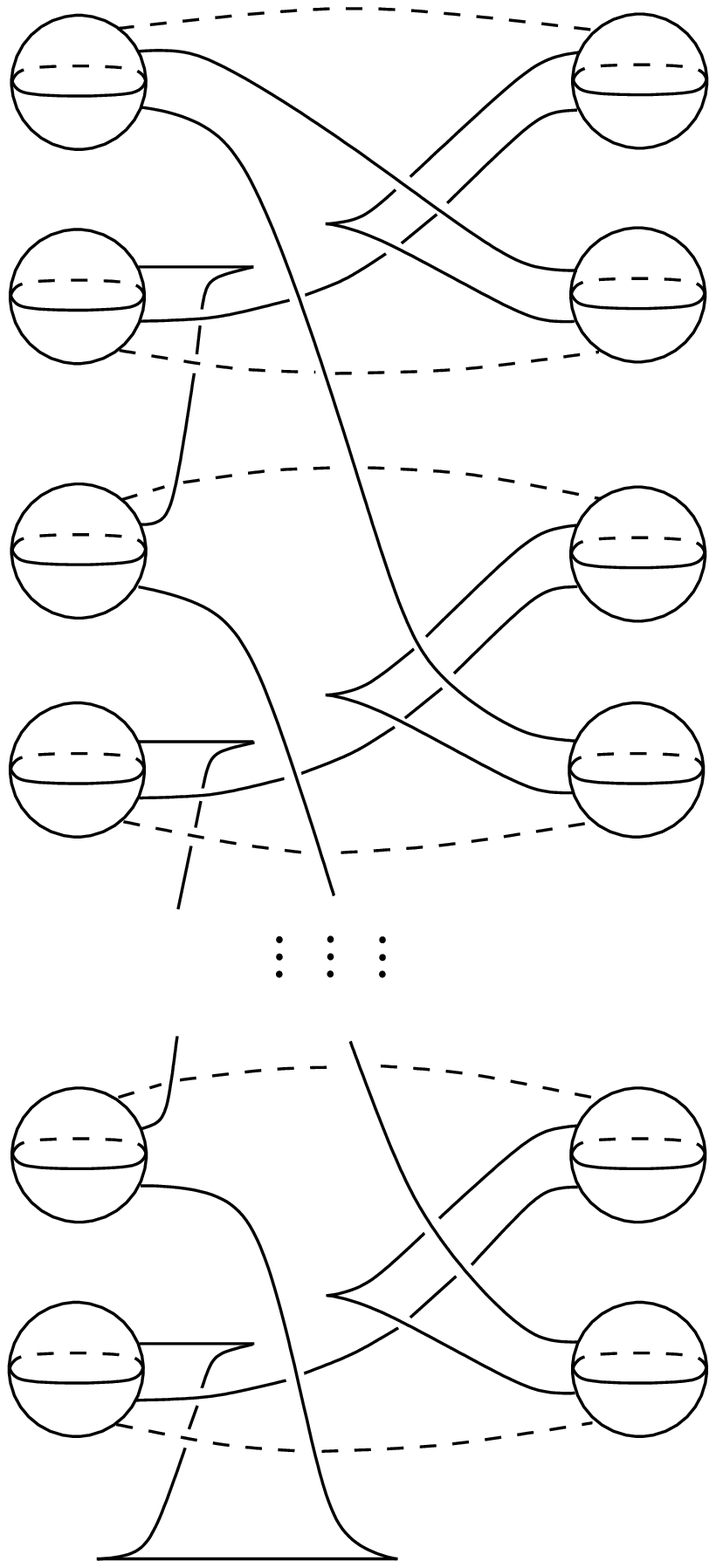}
  \captionof{figure}{Stein handlebody diagram for $\D T^*\S_g$ with $2g$ one-handles}
 \label{fig: orientable}
\end{minipage}%
\begin{minipage}{.5\textwidth}
  \centering
  \includegraphics[width=.7\linewidth]{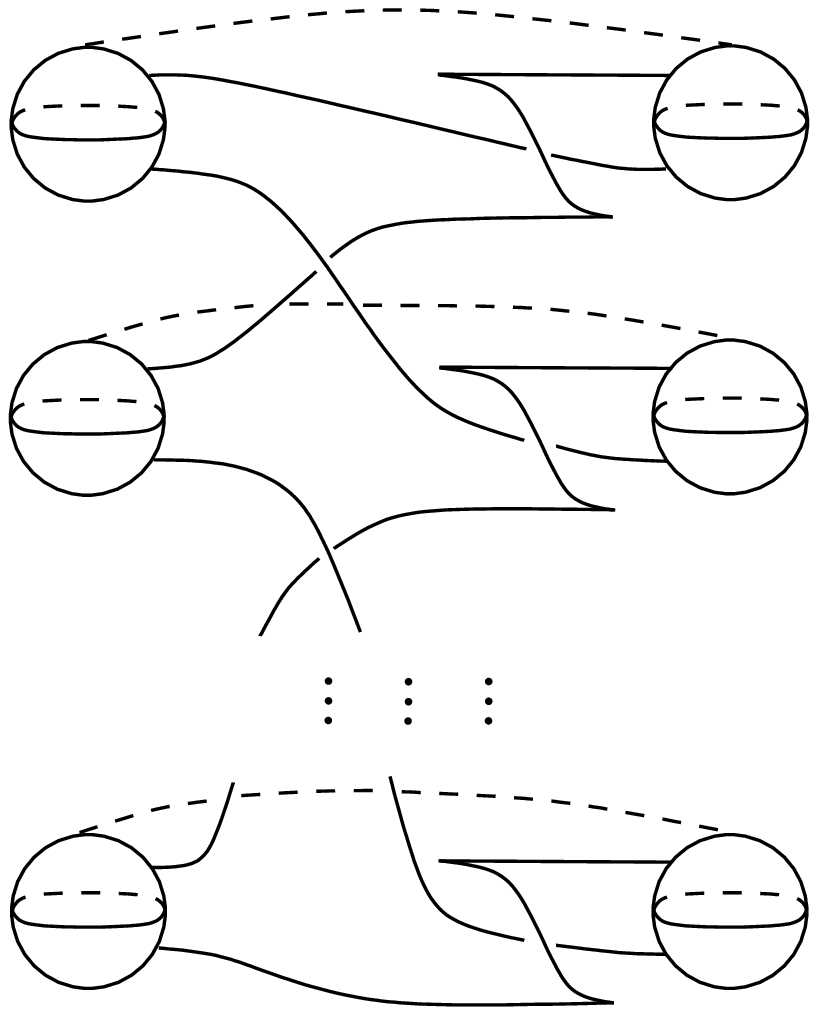}
  \captionof{figure}{Stein handlebody diagram for $\D T^*N_k$ with $k$ one-handles}
  \label{fig: nonorientable}
\end{minipage}
\end{figure}

\noindent (a) Suppose that $\S_g$ is a closed, connected, smooth and orientable surface of genus $g \geq1$. Then the Stein handlebody diagram depicted in Figure~\ref{fig: orientable} is symplectomorphic to $(\D T^*\S_g, \omega_{can})$ and its boundary is contactomorphic to $(\Sp T^*\S_g, \xi_{can})$.



(b) Suppose that $N_k$ is a closed, connected, smooth and nonorientable surface of genus $k \geq1$, i.e., $N_k = \#_k \mathbb{RP}^2$.  Then the Stein handlebody diagram depicted in Figure~\ref{fig: nonorientable} is symplectomorphic  to $(\D T^*N_k, \omega_{can})$ and its boundary is contactomorphic to  $(\Sp T^*N_k, \xi_{can})$.


\end{theorem}

\begin{remark} \label{rem: cases}

(i) The facts that the Stein handlebody diagram depicted in Figure~\ref{fig: orientable} is diffeomorphic to $\D T^*\S_g$ and the Stein handlebody diagram depicted in Figure~\ref{fig: nonorientable} is diffeomorphic to $\D T^*N_k$ were already proven by Gompf \cite{go}.

(ii) The  unit cotangent bundle $\Sp T^*\Sp^2$ is diffeomorphic to the real projective space $\mathbb{RP}^3$, and $\xi_{can}$ is the  unique tight contact structure on $\mathbb{RP}^3$, up to isotopy (cf. \cite{h}).  Moreover, McDuff \cite{mc} showed that any minimal symplectic filling of
        $(\mathbb{RP}^3, \xi_{can})$ is  diffeomorphic to $\D T^*\Sp^2$ and Hind \cite{hi} showed that $\D T^*\Sp^2$
is the unique Stein filling up to Stein homotopy. Furthermore, Wendl  \cite[Corollary 9.44]{w1} (based on his earlier work \cite{w}) showed that  any minimal strong symplectic
        filling of $(\mathbb{RP}^3, \xi_{can})$ is symplectic deformation equivalent
        to $(\D T^*\Sp^2, \omega_{can})$. A Stein structure on $\D T^*(\Sp ^2)$, which is diffeomorphic to the disk bundle over the $2$-sphere with Euler number $-2$, can be described by a single Stein handle attachment along a trivial Legendrian knot in the standard contact $3$-sphere. The boundary of this Stein domain is indeed contactomorphic to $(\Sp T^*(\Sp ^2), \xi_{can})$.

(iii) The unit cotangent bundle $\Sp T^* \T ^2$ is diffeomorphic to the
        $3$-torus $\T ^3$ and Eliashberg \cite{el} showed that $\xi_{can}$
        is the  unique strongly symplectically  fillable contact structure on $\T ^3$, up to contactomorphism.
        Moreover, according to Wendl \cite{w}, any minimal strong symplectic
        filling of $(\T ^3, \xi_{can})$ is symplectic deformation equivalent
        to $(\D T^* \T ^2\cong \T ^2 \times {\D}^2, \omega_{can})$.

(iv)  The unit cotangent bundle $\Sp T^*\mathbb{RP}^2 $ is diffeomorphic to the
        lens space  $L(4,1)$ and $\xi_{can}$
        is the  unique universally tight contact structure in $L(4,1)$, up to
        contactomorphism.  Moreover, McDuff
        \cite{mc} showed that $(L(4,1), \xi_{can})$ has two minimal symplectic fillings
        up to diffeomorphism: the disk cotangent bundle $\D T^*\mathbb{RP}^2$, which is a rational homology $4$-ball
        and  the disk bundle over the sphere with Euler number $-4$. Furthermore, Hind \cite{hi} showed that the uniqueness of the Stein fillings in each diffeomorphism class, up to Stein homotopy. More recently, Plamenevskaya and Van Horn-Morris \cite{pm} showed that, in each diffeomorphism class, there is a unique  minimal symplectic filling, up to symplectic deformation, based on the work of Wendl \cite{w}.

\end{remark}

\section{Upgrading the disk cotangent bundle to a Weinstein, and hence Stein filling}\label{sec: upgrade}

Let $q_1, q_2$ denote local coordinates on $S$, and $p_1, p_2$ denote dual coordinates for the cotangent fibers.
Then we have  $\lambda_{can}= \S \: p_idq_i$, and $\omega_{can}=d\lambda_{can}= \S \; d p_i \wedge d q_i$. It follows that $\Sigma\; p_i \partial_{p_i}$ is a Liouville vector field for the symplectic manifold $(\D T^*S, \omega_{can})$ transversely pointing out of $\partial (\D T^*S)$, which shows that $(\D T^*S, \omega_{can})$   is an exact  symplectic
        filling of its contact boundary  $(\Sp T^*S, \xi_{can}) $.  Now we briefly explain how this exact symplectic filling can be upgraded to the canonical Weinstein filling of $(\Sp T^*S, \xi_{can}) $ as described in  \cite[Example 11.12 (2)]{ce}. Fix any Riemannian metric on $S$ and a Morse function $f: S \to \R$. Let $X = \S \; p_i \partial_{p_i} + X_F$, where $X_F$ is the Hamiltonian vector field of $F = \lambda_{can}(\nabla f)$. Then, provided that $f$ is small enough,  $X$ is Liouville for $\omega_{can}$ and gradient-like for the Morse function $\phi(v)= 1/2 \| v\|_{\pi(v)} + f\circ \pi(v)$, where $\pi$ denotes the bundle projection $\D T^*S \to S$.  Thus, $(\D T^*S, \omega_{can}, X, \phi)$   is a Weinstein
        filling of  $(\Sp T^*S, \xi_{can}) $. Therefore, according to  \cite[Theorem 13.5]{ce}, $\D T^*S$ admits  a Stein domain structure $(J, \phi)$ such that the Weinstein domain associated to $(\D T^*S, J, \phi)$ is homotopic to $(\D T^*S, \omega_{can}, X, \phi)$. The main goal of this paper is to show that such a Stein domain structure on $\D T^*S$ is given by the handlebody diagrams in Figures~\ref{fig: orientable} and ~\ref{fig: nonorientable}, up to isotopy of $\D T^*S$.


\section{Weinstein homotopies}

\subsection{Orientable case} In this section we give a proof of Theorem~\ref{thm: gene} (a).

\begin{proof}

Let $(J_0, \phi_0)$ denote the Stein structure on $\D T^*\S_g$  given by the diagram in Figure~\ref{fig: orientable}.  We first observe that   $c_1 (  \D T^*\S_g , J_0) =0$ since the rotation number of the $2$-handle is zero (cf. \cite[Proposition 2.3]{go}).  On the other hand, if $J_{can}$ denotes an almost complex structure on $\D T^*\S_g$ which is  compatible  with  $\omega_{can}$, then $c_1 (  \D T^*\S_g , J_{can}) =0$ as well. In fact, any exact filling of $(\Sp T^*\S_g , \xi_{can}) $ has vanishing first Chern class \cite[Theorem 1.4]{lmy}.  Note that $J_{can}$ belongs to the unique homotopy class of almost complex structures on  $\D T^*\S_g$ compatible  with  $\omega_{can}$. Moreover,  since $H^2 (\D T^*\S_g; \Z) \cong \Z$ has no $2$-torsion, the homotopy class of  an almost complex structure on $\D T^*\S_g$  is determined by its first Chern class (see, for example, \cite[page 437]{gs}).    We conclude that the integrable almost complex structure $J_0$ is homotopic to the $\omega_{can}$-compatible almost complex structure $J_{can}$ on $\D T^*\S_g$. Now, by taking $V= \D T^*\S_g$, $\omega= \omega_{can}$, $X$ and $\phi$ as defined in Section~\ref{sec: upgrade},  and $J=J_0$, we deduce by  Theorem~\ref{thm: homoto} below that the Stein handlebody diagram depicted in Figure~\ref{fig: orientable} is Weinstein homotopic the canonical one, up to isotopy of $\D T^*\S_g$.

\begin{theorem}  \cite[Theorem 13.8]{ce} \label{thm: homoto} Let $(V, \omega, X, \phi)$ be a Weinstein manifold. Let $J$ be an integrable complex structure on $V$ which is homotopic to an almost complex structure compatible with $\omega$. Then there exists a diffeomorphism  $h: V \to V $ isotopic to the identity such that the function $\phi \circ h$ is $J$-convex and the Weinstein structure associated to the Stein structure $(h^*J, \phi)$ is homotopic to $(V, \omega, X, \phi)$, with fixed function $\phi$.   \end{theorem}

Since Weinstein homotopic manifolds are symplectomorphic  \cite[Corollary 11.21]{ce}, it follows that the Stein handlebody diagram depicted in Figure~\ref{fig: orientable} is indeed symplectomorphic to $(\D T^*\S_g, \omega_{can})$. Therefore,
the boundary the Stein handlebody diagram in Figure~\ref{fig: orientable} is contactomorphic to $(\Sp T^*\S_g, \xi_{can})$.\end{proof}

\begin{remark}  \label{rem: alt} Let $W_g$ be the Stein domain with boundary described by the Legendrian  handlebody diagram depicted in Figure~\ref{fig: orientable} and let $\xi_g$ denote the contact structure  induced on $\partial W_g$. An independent $3$-dimensional proof of the fact that $(\partial W_g, \xi_g)$ is contactomorphic to $(\Sp T^*\S_g, \xi_{can})$ can be given as follows.   Note that  $\partial W_g $ is the  circle bundle over $\S_g$ with Euler number $2g-2$, which is diffeomorphic to the unit cotangent bundle $\Sp T^*\S_g$. Let $\pi_g$ denote this circle fibration $\partial W_g \to \S_g$.   As shown by Lisca and Stipsicz \cite[Lemma 2.1]{ls}, $\xi_g$ has negative twisting number, and  thus  it is horizontal, i.e., $\xi_g$  is isotopic to a contact structure transverse to the circle fibers, by Honda's classification \cite[Theorem 2.11]{h1}  of tight contact structures on circle bundles over surfaces.  It follows that $\xi_g$  is universally tight by the work of Giroux \cite[Proposition 2.4 (c)]{g} and Honda \cite[Lemma 3.9]{h1}. As a matter of fact,  the twisting  number of $\xi_g$ is equal to $-1$ since there is a Legendrian knot $L$ in $(\partial W_g, \xi_g)$ as depicted in Figure~\ref{fig: legen}, which is isotopic to a fiber of $\pi_g$ (see  \cite[pages 289-290]{ls}). On the other hand, $\xi_{can}$ on $\Sp T^*\S_g$ is tangent to the fibers of the natural circle fibration $\Sp T^*\S_g \to \S_g$, by definition. Nevertheless,  it can be made horizontal by an arbitrarily small isotopy \cite[Proposition 1.4]{g} and  thus $\xi_{can}$ is also universally tight. Moreover, $\xi_{can}$ has twisting number  $-1$, by   \cite[Lemma 3.6]{g}. Since $\Sp T^*\S_g$ is diffeomorphic to $\partial W_g$ and there is a unique isomorphism class of universally tight contact structures on $\partial W_g$ with twisting number $-1$ (cf. \cite[Theorem 3.1 (c)]{g}), we conclude that $(\partial W_g, \xi_g)$ is contactomorphic to $(\Sp T^*\S_g, \xi_{can})$ for any $g \geq 1$.

\begin{figure}[h]
\centering
\relabelbox \small {\epsfxsize=1.5in
  \centerline{\epsfbox{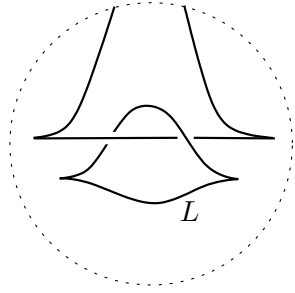}}}
\relabel{a}{{$L$}}

  \endrelabelbox
                        \caption{The Legendrian knot $L$ linking the $2$-handle at the bottom of Figure~\ref{fig: orientable} }
                        \label{fig: legen}
 \end{figure}

\end{remark}


\subsection{Nonorientable case} In this section we give the proof of Theorem~\ref{thm: gene} (b).

\begin{proof}  Let $V_k$ be the Stein domain with boundary described by the Legendrian handlebody diagram depicted in Figure~\ref{fig: nonorientable}. As observed in \cite[page 680]{go},  one can check via Kirby calculus that $V_k$ is diffeomorphic to the disk bundle over $N_k$  with Euler number $k-2 =-\chi(N_k)$, which is the disk cotangent bundle $\D T^*N_k$.    It follows that $\partial V_k$  is the circle bundle over $N_k$ with Euler number $k-2$, which is indeed diffeomorphic to $\Sp T^*N_k$.

Let $\xi_k$ denote the contact structure on $\partial N_k$ induced by the Stein handlebody diagram in Figure~\ref{fig: nonorientable}.  We first observe that $\xi_k$ has twisting number is $-1$, which follows from the fact that there is a Legendrian unknot $K$ in $(\partial N_k, \xi_k)$  which is isotopic to a fiber of the circle fibration of $\partial V_k$ over $N_k$, just as in the orientable case discussed above.  Note that there is a double cover  $ \S_{k-1} \to N_k$ and hence by pulling back the circle fibration and the contact structure we have a contact double cover $(M_k, \tilde{\xi}_k) \to (\partial V_k, \xi_k)$, where $M_k$ is the circle bundle over $\S_{k-1}$ with Euler number $2k-4$. The twisting number of $\tilde{\xi}_k$ is also $-1$ and thus  $(M_k, \tilde{\xi}_k)$ is universally tight just as in
Remark~\ref{rem: alt}. Therefore, we  conclude that $(\partial V_k, \xi_k)$ is universally tight as well, since it has a contact double cover which is universally tight.

On the other hand, $(\Sp T^*\S_{k-1}, \xi_{can})$ is the contact double cover of $(\Sp T^*N_k, \xi_{can})$, which implies in particular that $(\Sp T^*N_k, \xi_{can})$ is universally tight. Moreover, according to  \cite[Lemma 3.6]{g}, the twisting number of $(\Sp T^*N_k, \xi_{can})$ is $-1$, just as in the orientable case. Now we simply observe that the contact double cover $(\Sp T^*\S_{k-1}, \xi_{can})$ of $(\Sp T^*N_k, \xi_{can})$ is contactomorphic to the contact double cover $(M_k, \tilde{\xi}_k)$ of $(\partial V_k, \xi_k)$ by the proof given in Remark~\ref{rem: alt}.  Moreover, we may assume that this contactomorphism respects the circle fibrations (cf. \cite[Section 3 E] {g}) and hence yields a contactomorphism of  $(\Sp T^*N_k, \xi_{can})$ and  $(\partial V_k, \xi_k)$ after taking the quotients. This finishes the proof of our claim that $(\partial V_k, \xi_k) $ is contactomorphic to  $(\Sp T^*N_k, \xi_{can})$.

In the following, we prove the first assertion in Theorem~\ref{thm: gene} (b). Let $J_0$ denote the Stein structure on $\D T^*N_k$  given by the diagram in Figure~\ref{fig: nonorientable}, and let $J_{can}$ denote an almost complex structure compatible with $\omega_{can}$. Note that there is a unique homotopy class of $\omega_{can}$-compatible  almost complex structures and we claim that $J_0$ belongs to that class. Recall that almost complex structures on $\D T^*N_k$ correspond bijectively to $Spin^c$ structures. We denote by $\sp_{J}$ the $Spin^c$ structure on $\D T^*N_k$ associated to an almost complex structure $J$ (see, for example, \cite[Chapter 6]{ozst}).
Note that there is an injective map $$ Spin^c(\D T^*N_k) \stackrel{\psi}{\longrightarrow} Spin^c(\Sp T^*N_k)$$ such that for any almost complex structure $J$ on $\D T^*N_k$, we have  $\psi(\sp_J) = \t_{\xi_J}$, where  $\xi_J$ is the oriented $2$-plane field on $\Sp T^*N_k$ induced by $J$, and  $\t_{\xi_J}$ is the $Spin^c$ structure associated to $\xi_J$. By the second assertion in Theorem~\ref{thm: gene} (b), which we proved above, we have $\t_{\xi_{J_0}} = \t_{\xi_{J_{can}}}$, and thus  $\sp_{J_0} = \sp_{J_{can}}$, since $\psi$ is injective. This implies that $J_0$ is homotopic to $J_{can}$, proving our claim. Now, the first assertion in Theorem~\ref{thm: gene} (b) follows by  Theorem~\ref{thm: homoto} just as in the orientable case we discussed above.  \end{proof}



\section{Surgery diagrams for the canonical contact structures}

\begin{figure}[h]
\centering
\relabelbox \small {\epsfxsize=2.4in
  \centerline{\epsfbox{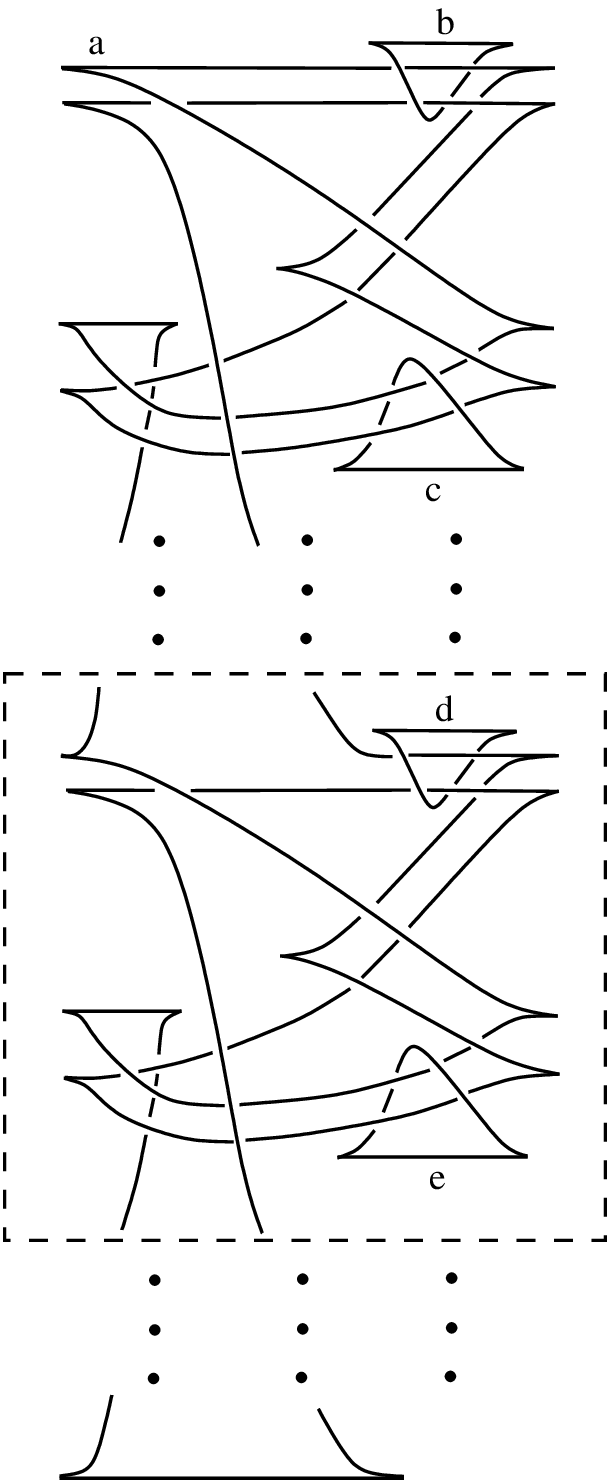}}}
\relabel{a}{{$-1$}}
  \relabel{b}{{$+1$}}
\relabel{c}{{$+1$}}
\relabel{d}{{$+1$}}
\relabel{e}{{$+1$}}
  \endrelabelbox
                        \caption{Surgery diagram for the canonical contact structure $\xi_{can}$ on the unit cotangent bundle of $\S_g$, where the dashed box is repeated $g-1$ times, for $g \geq 1$}
                        \label{fig: genusg}
 \end{figure}

\begin{figure}[h]
\centering
\relabelbox \small {\epsfxsize=2.7in
  \centerline{\epsfbox{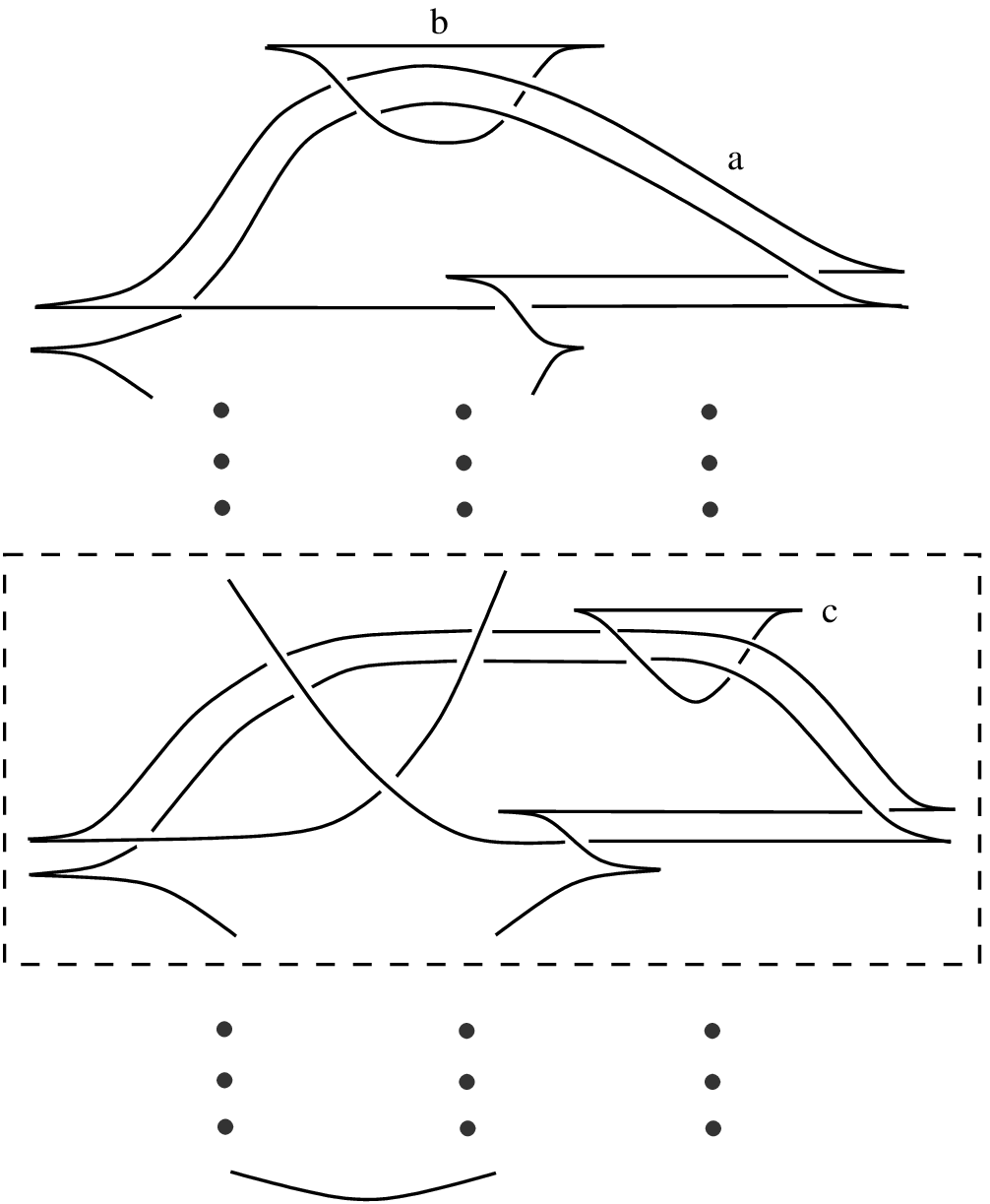}}}
\relabel{a}{{$-1$}}
  \relabel{b}{{$+1$}}
\relabel{c}{{$+1$}}

  \endrelabelbox
                        \caption{Surgery diagram for the canonical contact structure $\xi_{can}$ on the unit cotangent bundle of $N_k$, where the dashed box is repeated $k-1$ times, for $k \geq 1$}
                        \label{fig: genusk}
 \end{figure}

As explained in \cite[Theorem 5]{dg}, each $1$-handle in  a Stein handlebody diagram can be replaced by a contact $(+1)$-surgery along a Legendrian unknot, to obtain a surgery diagram of the contact boundary of the Stein handlebody. Therefore, as an immediate  corollary to Theorem 1.1, we obtain a surgery diagram for the contact $3$-manifold $(\Sp T^*S, \xi_{can})$ for any closed surface $S$. In Figure~\ref{fig: genusg}, we depicted the contact surgery diagram for $\xi_{can}$ on $\Sp T^*\S_g$ for $g \geq 1$. For $g=0$, see Remark~\ref{rem: cases} (ii). In Figure~\ref{fig: genusk},  we depicted the contact surgery diagram for $\xi_{can}$ on $\Sp T^* N_k$ for $k \geq 1$.

 \begin{figure}[h]
\centering
\relabelbox \small {\epsfxsize=1in
  \centerline{\epsfbox{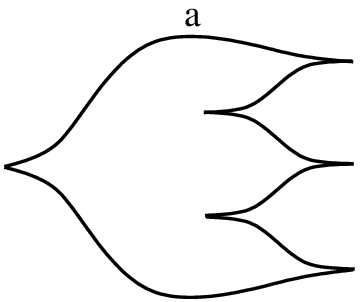}}}
\relabel{a}{{$-1$}}

  \endrelabelbox
                        \caption{A simple surgery diagram for $(\Sp T^*(\mathbb{RP}^2), \xi_{can})$}
                        \label{fig: rp}
 \end{figure}

\begin{remark} In Figure~\ref{fig: rp},  we depicted a simple surgery diagram for $(\Sp T^*(\mathbb{RP}^2), \xi_{can})$, without any contact $(+1)$-surgeries. The contact $3$-manifold described by the surgery diagram in Figure~\ref{fig: genusk} for $k=1$ is contactomorphic to the one in Figure~\ref{fig: rp}. As smooth $3$-manifolds the diffeomorphism between them can be given by a sequence of Kirby moves. Namely, after converting the diagram in Figure~\ref{fig: genusk} for $k=1$ into a smooth surgery diagram, we just perform a Rolfsen twist along the $0$-framed unknot and then blow-down the resulting $+1$-framed unknot to obtain a $-4$-framed unknot, which is the usual smooth diagram of $L(4,1)$. It is plausible that the aforementioned contactomorphism can be shown directly by using the set of moves in contact surgery diagrams introduced by Ding and Geiges \cite{dg}.

\end{remark}

{\bf {Acknowledgement}}: We would like to thank R. E. Gompf and A. I. Stipsicz for helpful comments on a draft of this paper.


\end{document}